\documentclass[
a4paper, 
oneside, 
]{amsart}

\usepackage[main=english, ngerman]{babel}
\usepackage[T1]{fontenc}
\usepackage[utf8]{inputenc}
\usepackage{csquotes}
\usepackage{lmodern}
\usepackage{cleveref}
\usepackage{stmaryrd}
\usepackage{graphicx}
\usepackage{mathrsfs} 
\usepackage{amsfonts}
\usepackage{color}
\usepackage{amsthm}
\usepackage{amsmath}
\usepackage{enumerate}
\usepackage[shortlabels]{enumitem}
\setlist[enumerate]{font=\normalfont}
\setlist[itemize]{font=\normalfont}
\usepackage{amssymb}
\usepackage{dsfont}
\usepackage{mathtools}
\usepackage{lineno}
\theoremstyle{break}
\newtheorem{thm}{Theorem}[section]
\newtheorem{cor}[thm]{Corollary}
\theoremstyle{definition}
\newtheorem{defi}[thm]{Definition}

\newtheorem{lem}[thm]{Lemma}
\newtheorem{rem}[thm]{Remark}
\newtheorem{prop}[thm]{Proposition}

\newcommand{\R}{\mathbb{R}}

\newcommand{\Q}{\mathbb{Q}}

\newcommand{\LL}{\mathrm{L}}
\newcommand{\X}{\mathrm{X}}
\newcommand{\Y}{\mathrm{Y}}
\newcommand{\der}{\mathrm{d}}
\newcommand{\dmu}{\,\mathrm{d}\mu}
\newcommand{\dnu}{\,\mathrm{d}\nu}
\newcommand{\dr}{\,\mathrm{d}r}
\newcommand{\ds}{\,\mathrm{d}s}
\renewcommand{\subset}{\subseteq}
\renewcommand{\phi}{\varphi}

\DeclareMathOperator{\id}{id}
\DeclareMathOperator{\Real}{Re}

\newcommand{\1}{\mathds{1}}

\newcommand{\fix}{\operatorname{fix}}
\newcommand{\argument}{\,\cdot\,}
\DeclarePairedDelimiter\norm{\lVert}{\rVert}
\DeclarePairedDelimiter\abs{\lvert}{\rvert}
\DeclareMathOperator{\lin}{span}
\DeclareMathOperator{\sign}{sign}

\usepackage[
    backend=biber,
    natbib=true,
    citestyle=numeric,
    url=false, 
    doi=false,
    isbn=false,
    eprint=false,
    date=year,
    maxbibnames=99
]{biblatex}

\DeclareFieldFormat[article]{number}{no. #1}
\renewbibmacro*{volume+number+eid}{%
  \printfield{volume}%
  \setunit*{\addnbspace}% NEW (optional); there's also \addnbthinspace
  \printfield{number}%
  \setunit{\addcomma\space}%
  \printfield{eid}}
\renewbibmacro{in:}{}

\title{Measure-preserving semiflows and one-parameter Koopman semigroups}
\author{Nikolai Edeko, Moritz Gerlach and Viktoria K\"uhner}
\address{Nikolai Edeko, Mathematisches Institut, Universit\"at T\"ubingen,
Auf der Morgenstelle 10, D-72076 T\"ubingen, Germany}
\email{nied@fa.uni-tuebingen.de}
\address{Moritz Gerlach\\Institut f\"ur Mathematik\\Universit\"at Potsdam\\Karl-Liebknecht-Stra\ss{}e 24--25\\14476 Potsdam}
\email{moritz.gerlach@uni-potsdam.de}
\address{Viktoria K\"uhner, Mathematisches Institut, Universit\"at T\"ubingen,
Auf der Morgenstelle 10, D-72076 T\"ubingen, Germany}
\email{viku@fa.uni-tuebingen.de}

\begin{document}

\begin{abstract}
For a finite measure space $\X$, we characterize strongly continuous Markov lattice semigroups 
on $\LL^p(\X)$ by showing that their generator $A$ acts 
as a derivation on the dense subspace $D(A)\cap\LL^\infty(\X)$. We then use this to 
characterize Koopman semigroups on $\LL^p(\X)$ if $\X$ is a standard probability space.
In addition, we show that 
every measurable and  measure-preserving flow on a standard probability space is 
isomorphic to a continuous flow on a compact Borel probability space.
\end{abstract}

\dedicatory{Dedicated to Rainer Nagel on the occasion of his $75$th birthday.} 

\maketitle

%\section{Introduction}

In this article we address mainly the following two issues. 
First, we characterize strongly continuous Markov lattice semigroups $(T(t))_{t\geq 0}$ on 
$\LL^p(\X)$ by properties of their generators for a finite measure space 
$\X = (X, \Sigma, \mu)$. We will show 
that a strongly continuous semigroup on $\LL^p(\X)$ is a Markov lattice semigroup
if and only if its generator $A$ acts as a derivation on $D(A)\cap\LL^\infty(\X)$, $\1\in D(A)$ and
the semigroup is locally bounded on $\LL^\infty(\X)$.
Similar results have been
established by R.\ Nagel and R.\ Derndinger in \cite[Satz 2.5]{nade} for semigroups on 
$\mathrm{C}(K)$, see also \cite[Section B-II.3]{lecture}, and recently by T.\ ter Elst and 
M.\ Lema\'{n}czyk  in \cite{terelst} for unitary groups on $\LL^2(\X)$.

Second, we show that such semigroups are always 
similar to a semigroup of Koopman operators. More precisely, we construct a compact 
space $K$ and a Borel measure $\nu$ such that $\LL^1(X,\Sigma,\mu)$ is isometrically 
Banach lattice isomorphic to $\LL^1(K,\nu)$ and, via this isomorphism, the semigroup 
$(T(t))_{t\geq 0}$ is similar to a semigroup of Koopman operators on $\LL^1(K,\nu)$ 
induced by a continuous semiflow $(\varphi_t)_{t\geq 0}$ on $K$. 
Furthermore, in case that the space $\LL^1(X,\Sigma,\mu)$ is separable, we  
show that $K$ can be chosen to be metrizable. Similar results have been already 
obtained for strongly continuous representations of locally compact groups on 
$\mathrm{L}^p(\X)$ as bi-Markov embeddings, see \cite[Theorem 5.14]{deJeu2017}.

The article is organized as follows. In the second part of the introduction, we 
specify our notation and recall some basic facts we use throughout the article. In  
\Cref{sec:characterization}, we prove our main result, \Cref{first}, that 
characterizes strongly Markov semigroups of lattice homomorphisms by the condition that 
their generator acts as a derivation, followed by a version for semigroups that 
are not necessarily Markov. In \Cref{sec:koopman} we then use these results to 
characterize Koopman semigroups and in particular obtain \cite[Theorem 1.1]{terelst} 
in which 
ter Elst and Lema\'{n}czyk proved a corresponding result for unitary operator groups as \Cref{cor:group}.
In \Cref{sec:topologicalmodel},
we turn to the construction of topological models. Finally, in Section 
\ref{sec:ergodicflows}, we consider ergodic, measure-preserving flows and give a new proof for 
the fact that they contain at most countably many non-ergodic mappings, provided
that their induced group of Koopman operators is strongly continuous on $\mathrm{L}^2(\X)$.
This has previously been proven in \cite[Theorem 1]{elements} for $\R^k$-actions.

\medskip

Let us recall some concepts and fix the notation used in this article. 
For a measure space $\X = (X,\Sigma,\mu)$ and $1\leq p\leq \infty$ we denote by
\begin{align*}
  \LL^p(\X) \coloneqq \LL^p(\X;\mathds{C}) = \LL^p(\X;\R)\otimes \mathrm{i}\LL^p(\X;\R)
\end{align*}
the corresponding
complex $\LL^p$-space. This is a complex Banach lattice in the sense of \cite[Definition 7.2]{erg}.
The lattice operations in $\LL^p(\X;\R)$ are denoted by $\vee$ and $\wedge$ and we write
$f_+$ and $f_-$ for the positive part $f\vee 0$ and negative part $(-f)\vee 0$ of a 
function $f\in\LL^p(\X;\R)$. We call $\X$ \emph{separable} if $\mathrm{L}^p(\X)$ is separable
for one (and hence all) $p\in[1,\infty)$.
When considering more than one measure on $X$, we may distinguish between a measurable function 
$f\colon X\to\mathbb{C}$ and its equivalence class with respect to a measure $\mu$ on $X$ by writing $[f]_\mu$
for the equivalence class of $f$. If $f\in\mathrm{L}^\infty(\X)$ is an essentially bounded function,
$\mathrm{M}_{f}$ will denote its associated multiplication operator on $\mathrm{L}^p(\X)$, $p\in [1,\infty]$.
Let $\X = (X, \Sigma, \mu)$ and $\Y = (Y, \Sigma', \mu')$ be finite measure spaces and $T\colon \LL^p(\X) \to \LL^p(\Y)$ a 
linear operator. The operator $T$ is called \emph{positive} if it is real, i.e., if 
$T\LL^p(\X;\R)\subseteq \LL^p(\Y;\R)$, and its restriction to the Banach lattice 
$\LL^p(\X;\R)$ is positive. The operator $T$ is said to be a \emph{lattice homomorphism} if $\abs{Tf} = T\abs{f}$ for each $f\in 
\LL^p(\X)$. 
In particular, every lattice homomorphism is positive and fulfills $T(f_+) = (Tf)_+$ and
$T(f_-) = (Tf)_-$ for each $f\in\mathrm{L}^p(\X;\R)$. The operator $T$ is called a \emph{Markov operator} 
if it is positive and $T\1_X = \1_Y$, and a \emph{bi-Markov} operator if, additionally, $T'\1_Y = \1_X$.
It is called a \emph{Koopman operator} if there is a measurable map $\varphi\colon Y\to X$
such that $\phi^{-1}$ maps null-sets into null-sets and $Tf = f\circ\varphi$ for all 
$f\in\LL^p(\X)$. In this case, we 
denote the operator by $T_\varphi$.
Note that every Koopman operator is a Markov lattice homomorphism.

Let $(T(t))_{t\geq 0}$ be a $C_0$-semigroup on $\mathrm{L}^p(\mathrm{X})$ for some  $1\leq p < \infty$. 
We denote its generator by $A$ and the domain of $A$ by $D(A)$.
The semigroup $(T(t))_{t\geq 0}$ is called a \emph{lattice semigroup} if each operator $T(t)$ is a lattice homomorphism.
If, additionally, for each $t\geq 0$ the operator $T(t)$ is a (bi-)Markov operator,
$(T(t))_{t\geq 0}$ is called a \emph{(bi-)Markov lattice semigroup}. It is called a \emph{Koopman 
semigroup}, if for each $t\geq 0$ the operator $T(t)$ is a Koopman operator.

Consider the equivalence relation 
\[
M\sim N \quad \text{if} \quad \mathds{1}_M = \mathds{1}_N \ \mu\text{-almost everywhere}
\]
on $\Sigma$. Then the set of equivalence classes $\Sigma(\X)\coloneqq\Sigma/{\sim}$ is called the
\emph{measure algebra} of the measure space $\X$ and is a Boolean algebra with respect to
the set operations union, intersection and complementation. 
For the sake of simplicity, we do not distinguish notationally between elements of $\Sigma$ and $\Sigma(\X)$.
A mapping
$\theta\colon \Sigma(\X) \to \Sigma(\X)$ is called a \emph{Boolean algebra homomorphism} if 
$\theta(\emptyset)=\emptyset$, $\theta(X)=X$ as well as
$\theta(A\cup B) = \theta(A)\cup \theta(B)$ and $\theta(A\cap B) = \theta(A)\cap \theta(B)$  for all $A,B\in \Sigma(\X)$.
If, in addition, a Boolean algebra homomorphism $\theta\colon \Sigma(\X)\to\Sigma(\X)$ satisfies $\mu(\theta(A))=\mu(A)$ for 
all $A\in \Sigma(\X)$,
then $\theta$ is called a \emph{measure algebra homomorphism}. Every measurable map $\varphi\colon X\to X$ 
such that $\phi^{-1}$ maps null-sets into null-sets induces a measure 
algebra homomorphism
$\varphi^*\colon\Sigma(\X)\to \Sigma(\X)$ via $A\mapsto \varphi^{-1}(A)$.
For further information on $\Sigma(\X)$ we refer to \cite[Section 6.1]{erg}. We call 
a measure space $(X, \Sigma, \mu)$ a \emph{Borel probability space} if $X$ can be equipped with 
a Polish topology such that $\Sigma$ is the corresponding Borel $\sigma$-algebra. 
A measure space $(X, \Sigma, \mu)$ is called a \emph{standard probability space} if there is a measurable, measure-preserving 
and essentially invertible map to a Borel probability space.

Finally, we call a linear operator $\delta$ on $\mathrm{L}^p(\mathrm{X})$ with domain 
$D(\delta)$ a \emph{derivation on $D(\delta)\cap\LL^\infty(\X)$} if
$D(\delta)\cap\LL^\infty(\X)$ is an algebra (with respect to the pointwise multiplication) 
and $\delta(f\cdot g)= \delta f \cdot g + f \cdot \delta g$ for all 
$f,g\in D(\delta)\cap\LL^\infty(\X)$.

\section{Characterization of Markov lattice semigroups on $\LL^p$-spaces}
\label{sec:characterization}

In this section we characterize strongly continuous Markov lattice semigroups on $\LL^p(\X)$-spaces,
where $\X = (X, \Sigma, \mu)$ is a finite measure space, by means of their generators. 
The following theorem is our main result.

\begin{thm}\label{first}
Let $A$ be the generator of a $C_0$-semigroup $(T(t))_{t\geq 0}$ on 
a space $\LL^p(\X)$, where $\X$ is a finite measure space
and $1\leq p < \infty$.
Then the following assertions are equivalent.
\begin{enumerate}[(i)]
\item $(T(t))_{t\geq 0}$ is a Markov lattice semigroup.
\item For every $t\geq 0$ there exists a Boolean algebra homomorphism 
$\theta_t \colon \Sigma(\X) \to \Sigma(\X)$ such that $T(t)\mathds{1}_M=\mathds{1}_{\theta_t(M)}$ 
for all $M\in \Sigma(\X)$.
\item The space $\LL^\infty(\X)$ is invariant under $(T(t))_{t\ge 0}$, the map
$t \mapsto \norm{T(t)}_{\mathscr{L}(\LL^\infty(\X))}$ is locally bounded, $\mathds{1}\in D(A)$  
and $A$ is a derivation on $D(A)\cap\LL^\infty(\X)$.
\end{enumerate}
\end{thm}

\begin{rem}
  \begin{enumerate}[(i)]
    \item   Given a finite measure space $\X$ and a bounded operator
            $S$ on $\LL^p(\X)$ such that $\LL^\infty(\X)$ is invariant 
            under $S$, it follows from the closed graph theorem that the restriction of 
            $S$ to $\LL^\infty(\X)$ is a bounded operator. Therefore, the map 
            $t \mapsto \norm{T(t)}_{\mathscr{L}(\LL^\infty(\X))}$ is well-defined in \normalfont{(iii)}.
            As will be shown in
            \Cref{lem:locbound}, the local boundedness condition is automatically fulfilled
            if $T$ is an operator group.
    \item A semigroup $(T(t))_{t\geq 0}$ satisfying \normalfont{(i)-(iii)} in \Cref{first}
    uniquely extends to a strongly continuous Markov lattice semigroup on $\mathrm{L}^q(\X)$ for
    each $1\leq q < \infty$ with 
    \begin{align*}
      \|T(t)\|_{\mathscr{L}(\mathrm{L}^q(\X))} = \|T(t)\|^{\frac{1}{q}}_{\mathscr{L}(\mathrm{L}^1(\X))},
    \end{align*}
    use \cite[Theorem 7.23]{erg}. In particular, it extends to the biggest $\mathrm{L}^q$-space, 
    $\mathrm{L}^1(\X)$, and we will therefore only consider semigroups on $\mathrm{L}^1(\X)$ in 
    \Cref{sec:topologicalmodel,sec:ergodicflows}. At this point, however, the dependence on 
    $p$ in \Cref{first} cannot be eliminated as easily since it is only clear a posteriori 
    that the implication (iii) $\implies$ (i) can be reduced to the case $p=1$.
  \end{enumerate}
\end{rem}

\medskip

As a preparation for the proof of \Cref{first}, recall the following lemma 
relating the algebra and the lattice structure of $\mathrm{L}^\infty(\X)$.
\begin{lem}\label{lem:multlatt}
  Let $\X$ be a finite measure space and 
  $T\colon\LL^\infty(\X)\to\LL^\infty(\X)$ be a 
  bounded linear operator satisfying $T\mathds{1} = \mathds{1}$. Then the follwing assertions are 
  equivalent.
  \begin{enumerate}[(i)]
    \item $T$ is multiplicative.
    \item $T$ is a $\mathrm{C}^*$-homomorphism.
    \item $T$ is a lattice homomorphism.
  \end{enumerate}
\end{lem}
\begin{proof}
  Obviously, (ii) implies (i).
  The equivalence of (ii) and (iii) can be found in \cite[Theorem 7.23]{erg}. (There, 
  the operator is assumed to be conjugation-preserving but this 
  assumption is trivially superfluous.)
  The implication (i) $\implies$ (ii) follows from \cite[Theorem 4.13]{erg}, the analogous 
  statement for spaces of continuous functions, by applying the Gelfand-Naimark theorem 
  \cite[Theorem 4.23]{erg}.
\end{proof}

The following continuity property will be essential for the proof of \Cref{first}.

\begin{lem}\label{fg}
Let $\X$ be a measure space and $B\subset \LL^\infty(\X)$ be bounded. Then the 
multiplication $\LL^p(\X)\times B \to \LL^p(\X)$, $(f, g) \mapsto fg$ is $\|\cdot\|_p$-continuous.
\end{lem}
\begin{proof}
Let $M$ be a bound for $B$. For $f, u\in\LL^p(\X)$, $g, v\in B$ and $c > 0$
\begin{align*}
  fg - uv 
  &= (f-u)g + u(g-v) \\
  &= (f-u)g + u\1_{[|u|\leq c]}(g- v) + u\1_{[|u|> c]}(g- v)
\end{align*}
and so 
\begin{align*}
  \limsup_{(f, g)\to (u, v)} \|fg - uv\|_p \leq 2M\left\|u\1_{[|u|> c]}\right\|_p \to 0 \quad (c\to \infty).
\end{align*}
\end{proof}

Next, we apply \Cref{fg} to a $C_0$-semigroup as in \Cref{first}.

\begin{cor}\label{productrule}
Let $A$ be the generator of a $C_0$-semigroup $(T(t))_{t\geq 0}$ on 
a space $\LL^p(\X)$, where $\X$ is a finite measure space
and $1\leq p < \infty$. In addition, suppose that the space 
$\LL^\infty(\X)$ is invariant under $(T(t))_{t\geq 0}$ and that the map 
$t \mapsto \norm{T(t)}_{\mathscr{L}(\LL^\infty(\X))}$ is locally bounded. Then, 
for all $f\in \LL^\infty(\X)$ and all $g\in \LL^p(\X)$ the function 
\begin{align*}
 [0,\infty) \to t\mapsto \mathrm{L}^p(\X), \quad t\mapsto T(t)f\cdot T(t)g
\end{align*}
is continuous. 
Moreover, for $f,g\in D(A)\cap \LL^\infty(\X)$ this function 
is differentiable 
and the product rule
\[
\frac{\der}{\der t} \left(T(t)f \cdot T(t)g\right) = T(t)Af \cdot T(t)g + T(t)f\cdot T(t)Ag
\]
holds.
\end{cor}
\begin{proof}
It suffices to prove the second part since the first is a consequence of \Cref{fg}.
Let $f,g \in D(A)\cap \LL^\infty(\X)$ and $t\geq 0$.
Use \Cref{fg} and differentiate to obtain
\begin{align*}
\frac{\der}{\der t} \left(T(t)f \cdot T(t)g\right)
&= \lim_{h\to 0} \frac{1}{h}  \bigl(\left[ T(t+h)f-T(t)f\right]T(t)g \\
&\qquad+ T(t+h)f\left[T(t+h)g-T(t)g\right]\bigr) \\
&= \lim_{h\to 0} \frac{1}{h} \bigl( \left[ T(t+h)f-T(t)f\right]\bigr)T(t)g \\
&\qquad+ \lim_{h\to 0}T(t+h)f\cdot \lim_{h\to 0} \frac{1}{h} \left[ T(t+h)g-T(t)g\right] \\
&= \left(\frac{\der}{\der t} T(t)f\right)\cdot T(t)g + T(t)f\cdot \left(\frac{\der}{\der t} T(t)g
\right)
\end{align*}
which proves the assertion.
\end{proof}

We are now able to prove \Cref{first}.

\begin{proof}[Proof of \Cref{first}.]
The equivalence (i) $\Leftrightarrow$ (ii) is proved almost exactly as in the case of 
bi-Markov operators, see \cite[Theorem 12.10]{erg}.

\medskip

To prove the implication (i) $\implies$ (iii), first note that
every operator $T(t)$ is positive. Since $T(t)\mathds{1}=\mathds{1}$ for each $t\geq 0$, this already 
implies that the semigroup preserves the subspace $\LL^\infty(\X)$ and the restriction 
of each $T(t)$ to $\LL^\infty(\X)$ is a contraction. In particular, it follows from 
\Cref{lem:multlatt} that every operator $T(t)$ is multiplicative on $\LL^\infty(\X)$.
By \Cref{productrule}, for every $f,g\in D(A)\cap \LL^\infty(\X)$ 
\begin{align*} \frac{\der}{\der t} T(t)(f\cdot g) &= \frac{\der}{\der t} \bigl( T(t)f \cdot T(t)g\bigr)\\
&= T(t)Af\cdot T(t)g + T(t)f\cdot T(t)Ag\\
&= T(t) \left[Af \cdot g + f \cdot Ag \right ] \, .
\end{align*}
In particular, this shows $f\cdot g \in D(A)$ and $A(f\cdot g) = Af\cdot g + f\cdot Ag$ for $t=0$.
This proves that $A$ is a derivation and clearly $\mathds{1}\in D(A)$. 
% The identity
% $A'\1 = 0$ is equivalent to $T(t)'\1 = \1$ for each $t \geq 0$ and hence follows because 
% $T(t)$ is bi-Markov.

\medskip

We now prove that (iii) implies (i). Because of the local boundedness of 
$t \mapsto \norm{T(t)}_{\mathscr{L}(\LL^\infty(\X))}$, there exists a constant $C>0$ 
such that
\[ \norm*{\frac{1}{t} \int_0^t T(s)f\ds}_\infty \leq C\norm{f}_\infty \]
for $0<t\leq 1$ and $f\in \LL^\infty(\X)$. This implies that 
$D\coloneqq D(A)\cap \LL^\infty(\X)$ is a dense subspace of  $\LL^p(\X)$.
We use this fact to show that each $T(t)$ is multiplicative on $\LL^\infty(\X)$.
For fixed $f,g\in D$ and $t>0$ consider the mapping
\[ s\mapsto \beta(s)\coloneqq T(t-s)[T(s)f\cdot T(s)g]\] 
on $[0,t]$.
Since $\beta(0)=T(t)(f\cdot g)$ and $\beta(t)=T(t)f\cdot T(t)g$, it suffices to show that $\beta$ is 
constant. To this end, consider the operator valued mappings 
$P,Q \colon [0,t]\to \mathscr{L}(\LL^p(\X))$ given by $P(s)=T(t-s)$ and 
$Q(s) = \mathrm{M}_{T(s)f} \circ  T(s)$, where $\mathrm{M}_{T(s)f}$ denotes the multiplication with the bounded 
function $T(s) f$. It follows from \Cref{productrule} that $Q$ is strongly continuous 
and that for each $h\in D$, $s\mapsto Q(s)h$ is differentiable with derivative
\[ \frac{\der}{\der s} Q(s)h = T(s)Af\cdot T(s)h + T(s)f\cdot T(s)Ah = A\bigl(T(s)f\cdot T(s)h\bigr)\, . \]
Here, the second equality follows from the fact that, by assumption,
$A$ is a derivation and $D$ is invariant under each $T(t)$. In particular, $D$ is invariant under $Q$.
Since $P$ is also strongly continuous and $s\mapsto P(s)h$ is differentiable for all $h\in D$,
it follows from \cite[Lemma B.16]{engnag} that
\[ \beta'(s) = - AT(t-s)[T(s)f\cdot T(s)g] + T(t-s)A[T(s)f\cdot T(s)g] = 0\]
for all $s\in [0,t]$. This shows that $\beta$ is constant and thus that every $T(t)$ is multiplicative 
on $D$.

\medskip

Since the multiplication with a fixed bounded function induces a bounded operator on $\LL^p(\X)$
and $D$ is $\|\cdot\|_p$-dense in $\LL^\infty(\X)$, fixing a function $g\in D$ and using a 
standard approximation argument shows that $T(f\cdot g) = T(t)f\cdot T(t)g$ for all 
$f\in \LL^\infty(\X)$ and $g\in D$. Fixing $f\in\LL^\infty(\X)$ and repeating the argument shows 
that $T(t)(f\cdot g)  = T(t)f\cdot T(t)g$ for all $f, g\in\LL^\infty(\X)$, so $T(t)$ is 
multiplicative on all of $\LL^\infty(\X)$.
Furthermore, $A\mathds{1}=0$ since $A$ is a derivation, i.e., $T(t)\mathds{1} = \mathds{1}$ for all $t\geq 0$.
Now \Cref{lem:multlatt} yields that every $T(t)$ is a lattice homomorphism on 
$\LL^\infty(\X)$ and hence, by density and continuity, also on $\LL^p(\X)$.
\end{proof}

\begin{rem}\label{oneremark}
 The assumption $\1\in D(A)$ in (iii) of \Cref{first} is automatically satisfied
 if $T(t)$ is an isometry for each $t \geq 0$.
 To see this, note that
 in the proof of \Cref{first} the assumption $\1\in D(A)$ was only important for the 
 implication (iii) $\implies$ (i).
 There, we showed that all the operators of the semigroup are multiplicative
 on $\LL^\infty(\X)$ and therefore map characteristic functions to 
 characteristic functions. If $T(t)$ is an isometry, it follows that $T(t)\1 = \1$ for each $t\geq 0$
 and hence $\1\in D(A)$. This assumption is also fulfilled if $T(t)'\1 = \1$ for each $t\geq 0$, 
 since then
  $\langle T(t)\1, \1\rangle = \langle \1, \1\rangle = \mu(X)$ and so,
 $T(t)\1$ being a characteristic function, $T(t)\1 = \1$ for each $t\geq 0$.
\end{rem}

As a corollary of \Cref{first}, we also obtain the following
characterization of bi-Markov lattice semigroups.

\begin{cor}\label{bimarkovchar}
  Let $A$ be the generator of a $C_0$-semigroup $(T(t))_{t\geq 0}$ on 
  a space $\LL^p(\X)$, where $\X$ is a finite measure space
  and $1\leq p < \infty$. Then the following assertions are equivalent.
  \begin{enumerate}[(i)]
    \item $(T(t))_{t\geq 0}$ is a bi-Markov lattice semigroup.
    \item For every $t\geq 0$ there exists a measure algebra homomorphism 
    $\theta_t \colon \Sigma(\X) \to \Sigma(\X)$ such that $T(t)\mathds{1}_M=\mathds{1}_{\theta_t(M)}$ 
    for all $M\in \Sigma(\X)$.
    \item The space $\LL^\infty(\X)$ is invariant under $(T(t))_{t\ge 0}$, the map
    $t \mapsto \norm{T(t)}_{\mathscr{L}(\LL^\infty(\X))}$ is locally bounded,
    $A$ is a derivation on $D(A)\cap\LL^\infty(\X)$ and $A'\1 = 0$.
  \end{enumerate}
\end{cor}
\begin{proof}
 For the equivalence of (i) and (ii), the reader is again referred to \cite[Theorem 12.10]{erg}.
 The equivalence
 of (i) and (iii) follows from \Cref{first} and \Cref{oneremark} since 
 $T(t)'\1 = \1$ for all $t\geq 0$ is equivalent to $A'\1 = 0$.
\end{proof}

In the following we discuss $C_0$-semigroups of lattice homomorphisms on $\LL^p(\X)$ that are not 
necessarily Markov. We show that their generator is a derivation perturbed by a bounded 
multiplication operator.
\begin{thm}\label{second}
Let $A$ be the generator of a $C_0$-semigroup $(T(t))_{t\geq 0}$ on 
a space $\LL^p(\X)$, where $\X$ is a finite measure space
and $1\leq p < \infty$. Assume that $\mathds{1}\in D(A)$ and 
$q\coloneqq A\mathds{1}\in \LL^\infty(\X)$. Then the following assertions are equivalent.
\begin{itemize}
\item[(i)] $(S(t))_{t\geq 0}$ is a lattice semigroup.
\item[(ii)] The function $q$ is real-valued and $A = B + q$ where $B$ is the 
generator of a Markov lattice $C_0$-semigroup $(T(t))_{t\geq 0}$ on 
$\LL^p(\X)$.
\end{itemize}
If {\normalfont{(ii)}} holds, then 
\begin{align}
\label{eqn:S(t)formula}
S(t)f = \exp{\left(\int_0^t \! T(s)q \ds\right)} \cdot T(t)f 
\end{align}
for all $t\geq 0$ and $f\in \LL^p(\X)$.
\end{thm}
\begin{proof} 
To show the equivalence of (i) and (ii), we first recall 
from \cite[Theorem III.1.3]{engnag} that $B\coloneqq A-q$ is a generator of a $C_0$-semigroup 
$(T(t))_{t\geq 0}$ on $\LL^p(\X)$
because $B$ is a bounded perturbation of $A$. Since $B \mathds{1} = 0$,  
$T(t)\mathds{1} = \mathds{1}$ for all $t\geq 0$.
Now it follows from \cite[Corollary C-II.5.8]{lecture} (Kato's identity) that $(S(t))_{t\geq 0}$ is 
a lattice semigroup if and only if $D(A)$ is a sublattice of $\LL^p(\X)$ and
\begin{align*}
  A\abs{f}=\Real(\sign(f)Af)
\end{align*}
for all $f\in D(A)$. Since (i) implies that $q$ is real-valued, $A$ satisfies 
this condition if and only if $B$ does, which proves the equivalence of the assertions 
(i) and (ii).

Now assume that (ii) holds.
% and consider the operators $S(t)\in \mathscr{L}(\LL^p(\X))$ given by
% \[S(t)f\coloneqq \exp{\left(\int_0^t \! T(s)q \ds\right)}T(t)f\]  for each $t\geq 0$.
Since each $T(t)$ is multiplicative on $\LL^\infty(\X)$ by \Cref{lem:multlatt}, 
$T(t)\exp(g) = T(t)\sum_{n=0}^{\infty}\frac{g^n}{n!}= \exp(T(t)g)$ for each $g\in \LL^\infty(\X)$.
Using this, one proves by induction that
\begin{align*}
  \left(T(t)\mathrm{e}^{t\mathrm{M}_q}\right)^n = \mathrm{M}_{\exp\left(t\sum_{j=1}^n T(jt)q\right)}T(nt)
\end{align*}
for each $n \geq 1$. Replace $t$ by $\frac{t}{n}$ and note that 
\begin{align*}
  \sum_{j=1}^n \frac{t}{n} T\left(\frac{jt}{n}\right)q 
  \xrightarrow[n\to\infty]{\|\cdot\|_p} \int_0^t T(s)q\ds.
\end{align*}
By \Cref{fg}, one obtains the convergence
\begin{align*}
  \mathrm{M}_{\exp\left(\sum_{j=1}^n \frac{t}{n} T\left(\frac{jt}{n}\right)q\right)}
  \xrightarrow[n\to\infty]{} \mathrm{M}_{\exp\left(\int_0^t T(s)q\ds\right)}
\end{align*}
of multiplication operators in the strong operator topology on 
$\mathscr{L}(\LL^p(\X))$. By the Trotter product formula
\begin{align*}
  S(t)f = \lim_{n\to\infty}\left[T\left(\frac{t}{n}\right)\exp\left(\frac{t}{n}\mathrm{M}_q\right)\right]^nf
  = \exp{\left(\int_0^t \! T(s)q \ds\right)} \cdot T(t)f 
\end{align*}
for all $f\in\LL^p(\X)$.
\end{proof}

\section{Koopman semigroups on $\LL^p$-spaces}\label{sec:koopman}
Every Koopman semigroup on an $\LL^p(\X)$-space is a Markov lattice semigroup but  
the converse is, in general, not true. However, it does hold if $\X$ is a standard 
probability space: In the case of bi-Markov lattice homomorphisms, this is a classical theorem 
by von Neumann, cf. \cite[Theorem 7.20]{erg}. Below, we give an operator-theoretic 
proof extending this theorem to Markov lattice homomorphisms on $\LL^p$-spaces. We then 
relate this to the results on semigroups from the previous section.

\begin{thm}\label{neumann}
 Let $\mathrm{X} = (X, \Sigma_X, \mu_X)$ and $\Y = (Y, \Sigma_Y, \mu_Y)$ be standard probability spaces
 and $T\colon \LL^p(\X) \to \LL^p(\Y)$, 
 $1\leq p \leq \infty$, a Markov lattice 
 homomorphism (not necessarily bi-Markov). Then there is a measurable map 
 $\varphi\colon Y\to X$ such that $T = T_\varphi$. If $\vartheta\colon Y\to X$ is another
 such map, then $\varphi = \vartheta$ $\mu_Y$-almost everyhwere.
\end{thm}
\begin{proof}
 We will factorize $T$ in order to apply von Neumann's theorem 
 for bi-Markov lattice homomorphisms. Since $\mathrm{L}^\infty(\X)$ is dense in $\mathrm{L}^p(\X)$
 for any $p\in [1, \infty]$, we only need to prove the case $p=\infty$. Moreover, 
 assume without loss of generality that $\X = (K, \mathcal{B}(K), \mu_K)$ and 
 $\Y = (L, \mathcal{B}(L), \mu_L)$ are Borel probability spaces. 
 
 \medskip
 
 Set $\nu \coloneqq (T'\1_Y)\mu$ and note that 
 $(K, \mathcal{B}(K), \nu)$ is again a Borel probability space. 
 Since $\nu \ll \mu$, the map 
 $P\colon \LL^p(K, \mathcal{B}(K), \mu_K) \to \LL^p(K, \mathcal{B}(K), \nu)$,
 $[f]_{\mu_K} \mapsto [f]_\nu$ is a well-defined, surjective and bounded operator.
 Moreover, for $f\in\mathrm{L}^\infty(K, \mathcal{B}(K), \mu_K)$, 
 \begin{align}\label{factorbimarkov}
   \int_K |f| \dnu = \int_L T|f| \dmu_L = \int_L |Tf|\dmu_L 
 \end{align}
 and so $\ker T \subset \ker P$. Therefore, there is an operator 
 $\hat{T}\colon\mathrm{L}^\infty(K, \mathcal{B}(K), \nu) \to \LL^\infty(L, \mathcal{B}(L), \mu_L)$
 such that $T = \hat{T}P$. Since $P$ is surjective, $\hat{T}$ also is
 a Markov lattice homomorphism and it follows from (\ref{factorbimarkov}) that 
 $\hat{T}$ is, in fact, bi-Markov.
 
 Now, 
 von Neumann's theorem shows that $\hat{T}[f]_\nu = [f\circ\varphi]_{\mu_L}$
 for a measurable and measure-preserving map 
 $\varphi\colon (L, \mathcal{B}(L), \mu_L)\to(K,\mathcal{B}(K),\nu)$
 and so $T = T_\phi$. The proof that $\phi$ is unique almost everyhwere is 
 the same as for measure-preserving maps, so we refer the reader to, e.g., 
 \cite[Lemma 6.9]{erg}.
\end{proof}

\begin{cor}\label{stand1}
Let $A$ be the generator of a $C_0$-semigroup $(T(t))_{t\geq 0}$ on 
a space $\LL^p(\X)$, where $\X = (X, \Sigma, \mu)$ is a standard probability space
and $1\leq p < \infty$. Then the equivalent assertions {\normalfont(i)}, 
{\normalfont(ii)} and {\normalfont(iii)} of \Cref{first} are also equivalent to 
\begin{itemize}
\item[(iv)] There exists a family $(\varphi_t)_{t\geq0}$ of measurable maps on $X$ 
such that $T(t)f=f\circ \varphi_t$ for all $f\in \LL^p(\X)$ and $t\geq 0$.
\end{itemize}
\end{cor}
\begin{proof}
If assertion (i) of \Cref{first} holds, we obtain assertion (iv) by \Cref{neumann} below.
Conversely, if (iv) holds, then every $T(t)$ is a Markov lattice homomorphism, thus 
assertion (i) holds.
\end{proof}

\begin{rem}
  Given a Koopman semigroup $(T_{\phi_{t}})_{t\geq 0}$ on $\mathrm{L}^p(\X)$ as in 
  \Cref{stand1}, it is immediate from the semigroup property and \Cref{neumann} 
  that $\phi_0 = \id_X$ $\mu$-almost everywhere and $\phi_t\circ\phi_s = \phi_{t+s}$ 
  $\mu$-almost everywhere for all $s, t\geq 0$. Therefore, the family 
  $(\phi_t)_{t\geq 0}$ forms a \emph{semiflow} modulo null-sets, see \Cref{sec:topologicalmodel}. 
  Note, however, that 
  it can in general not be made into a semiflow by simply discarding a null-set
  since the identity $\phi_t\circ\phi_s = \phi_{t+s}$ might hold outside of a null-set 
  depending on $s$ and $t$.
\end{rem}

\begin{cor}\label{bimarkovkoopman}
  Let $A$ be the generator of a $C_0$-semigroup $(T(t))_{t\geq 0}$ on 
a space $\LL^p(\X)$, where $\X = (X, \Sigma, \mu)$ is a standard probability space
and $1\leq p < \infty$. Then the equivalent assertions {\normalfont(i)}, 
{\normalfont(ii)} and {\normalfont(iii)} of \Cref{bimarkovchar} are also equivalent to 
  \begin{enumerate}
    \item[(iv)] There exists a family $(\varphi_t)_{t\geq0}$ of measurable and  
    measure-preserving maps on $X$ such that $T(t)f=f\circ \varphi_t$ for all 
    $f\in \LL^p(\X)$ and $t\geq 0$.
  \end{enumerate}
\end{cor}
\begin{proof}
 Assume (i) of \Cref{bimarkovchar}. Then \Cref{stand1} shows that there are measurable maps
 $\varphi_t\colon X\to X$ such that $T(t)f = f\circ\varphi_t$. Moreover, 
 for $M\in\Sigma$
  \begin{align*}
   \mu(\phi_t^{-1}(M)) 
   = \langle \mathds{1}_{\phi_t^{-1}(M)}, \1_X\rangle
   = \langle T(t)\mathds{1}_M, \1_X\rangle
   = \langle \mathds{1}_M, \1_X\rangle
   = \mu(M)
 \end{align*}
 and so each $\phi_t$ is 
 measure-preserving. On the other hand, (iv) implies (ii) with 
 $\theta_t = \varphi_t^*$.
\end{proof}

The next lemma shows that in the case of groups, the boundedness assumption
in \Cref{first} {\normalfont{(iii)}} is superfluous. This allows us to recover 
\cite[Theorem 1.1]{terelst} as \Cref{cor:group}.
The idea of the following proof is based on the proof of 
\cite[Theorem 2.5]{terelst}.

\begin{lem}\label{lem:locbound} Let $\X$ be a finite measure space 
  and $(T(t))_{t>0}$ be a 
  semigroup on $\LL^\infty(\X)$, strongly continuous with respect to $\|\cdot\|_p$ 
  where $1\leq p < \infty$. Then the mapping 
  $t \mapsto \norm{T(t)}_{\mathscr{L}(\LL^\infty(\X))}$ is locally bounded. 
\end{lem}
\begin{proof}
  For $f\in \LL^\infty(\X)$ and $t > 0$, setting $q \coloneqq 1 - \frac{1}{p}$
  \begin{align*}
    \norm{T(t)f}_\infty
    &= \sup\left\{ \abs{\langle T(t)f, g\rangle} \colon g\in \LL^1(\X) \text{ and } \norm{g}_1 = 1\right\} \\
    &= \sup\left\{ \abs{\langle T(t)f, g\rangle} \colon g\in \LL^q(\X) \text{ and } \norm{g}_1 = 1\right\} 
  \end{align*}
  Since $(T(t))_{t>0}$ is strongly continuous with respect to $\|\cdot\|_p$, $\langle T(\argument)f, g\rangle$
  is continuous for $g\in \LL^q(\X)$.
  Hence, $\norm{T(\argument)f}_\infty$ is lower semicontinuous, being the supremum of continuous
  functions. Therefore,  the sets 
  $A_n \coloneqq \{ t > 0 \colon \norm{ T(t)f}_\infty \leq n\}$ are closed.
  Since $\R_{>0} = \bigcup_n A_n$, Baire's category theorem yields that for 
  every $k\in\mathbb{N}$, there
  are $0< a_k < b_k \leq \frac{1}{k}$ and $M_k\in\mathbb{N}$ such that $[a_k, b_k]\subset A_{M_k}$. By
  the semigroup property, 
  \begin{align*}
    \sup_{s\in[a_k+t, b_k+t]}  \norm{T(s)f}_\infty 
    &= \sup_{s\in[a_k+t, b_k+t]} \norm{T(t)T(s-t)f}_\infty \\
    &\leq \sup_{s\in[a_k, b_k]} \norm{T(t)}_\infty\norm{T(s)f}_\infty \leq M_k\norm{T(t)}_\infty
  \end{align*}
  for each $t>0$.
  This shows that $\norm{T(\argument)f}_\infty$ is locally bounded for each $f\in\LL^\infty(\X)$
  and the claim now follows by the principle of uniform boundedness.
\end{proof}

%%%%%%%%%%%%%%%%%%%%%%% terelst
 
\begin{cor}[{\cite[Theorem 1.1]{terelst}}]\label{cor:group}
Let $A$ be the generator of a unitary C$_0$-group $(T(t))_{t\in\R}$ on 
$\LL^2(\X)$ where $\X = (X, \Sigma, \mu)$ is a standard probability space.
Then the following 
assertions are equivalent.
\begin{itemize}
\item[(i)]  For every $t\in\R$ there exists an essentially invertible measurable and measure-preserving 
map $\varphi_t \colon X \to X$ such that 
$T(t)f=f\circ\varphi_t$ for all $f\in\LL^2(\X)$.
\item[(ii)] 
The space $\LL^\infty(\X)$ is invariant under $(T(t))_{t\geq 0}$ and $A$ is a derivation
on $D(A)\cap\LL^\infty(\X)$.
\end{itemize}
\end{cor}
\begin{proof} 
The implication {\normalfont{(i)}} $\implies$ {\normalfont{(ii)}} is a consequence of 
\Cref{bimarkovchar}. 
In order to prove the converse implication, we observe that it follows from \Cref{lem:locbound}
and \Cref{oneremark} that 
$A$ and $(T(t))_{t\geq 0}$ as well as $-A$ and $(T(-t))_{t\geq 0}$ fulfill condition (iii) in \Cref{first}.
\Cref{stand1} therefore shows that $T(t) = T_{\varphi_t}$ for 
measurable maps $\varphi_t\colon X\to X$ and $t\geq 0$. The essential invertibility of the maps $\varphi_t$ follows from 
\cite[Proposition 7.12]{erg} and \cite[Corollary 7.21]{erg}. Also, since each $T(t)$ is unitary and
a Markov operator, one shows as in \Cref{bimarkovkoopman} that each $\phi_t$ is measure-preserving.
\end{proof}

\begin{cor}\label{stand2}
Let $A$ be the generator of a $C_0$-semigroup $(S(t))_{t\geq 0}$ on 
a space $\LL^p(\X)$, where $\X$ is a standard probability space
and $1\leq p < \infty$, such that $\1\in D(A)$ and 
$q \coloneqq A\mathds{1} \in \LL^\infty(\X)$. Then 
$(S(t))_{t\geq 0}$ is a lattice semigroup if and only if $q\in\LL^\infty(\X; \R)$ and 
there exists a family $(\varphi_t)_{t\geq 0}$ of measurable maps on $X$ corresponding 
to a strongly continuous Koopman semigroup on $\LL^p(\X)$ such that 
\begin{align}
\label{eqn:S(t)formulaflow}
S(t)f = \exp{\left(\int_0^t  q\circ \varphi_s \ds\right)} \cdot (f\circ \varphi_t)
\end{align}
for all $f\in \LL^p(\X)$ and $t\geq 0$.
\end{cor}
\begin{proof}
If $(S(t))_{t\geq 0}$ is a lattice semigroup, it follows from \Cref{second} that there exists a 
Markov lattice semigroup $(T(t))_{t\geq 0}$ on $\LL^p(\X)$ with generator $(A-q,D(A))$ such that 
\eqref{eqn:S(t)formula} holds. The representation \eqref{eqn:S(t)formulaflow} hence follows from 
\Cref{stand1}. Conversely, every semigroup of the form \eqref{eqn:S(t)formulaflow} with a real-valued
$q$ is a lattice semigroup.
\end{proof}

\section{Topological Model}
\label{sec:topologicalmodel}

We have seen in \Cref{stand1} that on a standard probability space $\X=(X,\Sigma,\mu)$ every 
strongly continuous Markov lattice semigroup $(T(t))_{t\geq 0}$ on $\LL^1(\X)$ is induced by a family 
$(\varphi_t)_{t\geq0}$ of measurable maps on $X$. Since $(T(t))_{t\geq 0}$ is a semigroup, one has 
$\varphi_0 = \mathrm{id}_X$ and $\varphi_s\circ\varphi_t = \varphi_{s+t}$ almost everywhere, 
using the uniqueness in \Cref{neumann}.
Recall
that a family $(\varphi_t)_{t\geq 0}$ of maps on a set $X$ is called a \emph{semiflow} if 
$\varphi_0 = \mathrm{id}_X$ and $\varphi_t\circ\varphi_s = \varphi_{s+t}$ for $s, t\geq 0$. 
A semiflow on a 
topological space $X$ is called \emph{continuous} if the map $\Phi\colon X\times\R_+\to X$, 
$(x, t)\mapsto \varphi_t(x)$ is 
continuous. Similarly, a semiflow on a measure space is called \emph{measurable} if 
$\Phi$ is measurable. It is called measure-preserving if for each 
$t\geq 0$ the map $\varphi_t$ is measurable and measure-preserving. In that case, there is an induced
semigroup of operators $(T_{\varphi_t})_{t\geq 0}$ on $\LL^1(\X)$ called the \emph{Koopman semigroup}
induced by the semiflow. This semigroup is weakly measurable if the semiflow is measurable but 
even then it need not be strongly continuous: Consider, e.g., the semiflow on 
$(\{0,1\}, \mathcal{P}(\{0,1\}), \frac{1}{2}(\delta_0 + \delta_1))$ with 
$\phi_0 = \id_{\{0,1\}}$ and $\phi_t \equiv 0$ for each $t>0$.
For measurable and measure-preserving \emph{flows} 
$(\varphi_t)_{t\in\R}$ (defined analogously to semiflows) on separable measure spaces, however, 
the induced Koopman group is in fact strongly continuous.

\begin{prop}\label{measflow}
  Let $\mathrm{X} = (X, \Sigma, \mu)$ be a separable, finite measure space and 
$(\varphi_t)_{t\geq 0}$ a measurable and measure-preserving semiflow on $X$.
    Then
    the induced Koopman semigroup $(T(t))_{t\geq 0}$ on $\mathrm{L}^p(\mathrm{X})$ is 
    strongly continuous on $(0,\infty)$ for $1 \leq p < \infty$.
    If $(\varphi_t)_{t\in\R}$ is a measurable and measure-preserving flow on $X$,
    the induced Koopman group on $\mathrm{L}^p(\mathrm{X})$ is 
    strongly continuous on all of $\R$ for $1 \leq p < \infty$.
\end{prop}
\begin{proof}
  Take $f\in\mathrm{L}^p(\mathrm{X})$. Then the function $f\circ\Phi$ is
  measurable on $X\times\R$ and for each $g\in\mathrm{L}^q(\mathrm{X})$ with 
  $\frac{1}{q} = 1 - \frac{1}{p}$ 
  the map
  \begin{align*}
     t \mapsto \langle T(t)f, g\rangle = \int_X f(\Phi(x, t))g(x)\dmu(x)
  \end{align*}
  is measurable (cf. \cite[Theorem 1.7.15]{TaoMeasure}). Therefore, the semigroup $(T(t))_{t\geq 0}$ 
  is weakly measurable and even strongly
  measurable by \cite[Theorem 3.5.5]{HillePhillips} since $\LL^p(\X)$ is separable. Since strongly 
  measurable semigroups are strongly continuous on $(0, \infty)$ by \cite[Theorem 10.2.3]{HillePhillips},  
  $(T(t))_{t\geq 0}$ is indeed strongly continuous. The case of flows follows immediately.
\end{proof}

We now show that for such measurable flows, one can construct a continuous flow on a 
compact metric space with an invariant probability measure such that the two flows are
isomorphic. This will be done by first proving that every strongly continuous Markov lattice semigroup is similar
to a Koopman semigroup induced by a continuous semiflow. An analogous 
result was recently proved for bi-Markov lattice embedding representations of locally compact
groups by de Jeu and Rozendaal, see \cite[Theorem 5.14]{deJeu2017}.

\begin{defi}
  Let $\X = (X, \Sigma, \mu)$ and $\Y = (Y, \Sigma', \mu')$ be finite measure
  spaces. We say that two measurable semiflows $(\varphi_t)_{t\geq 0}$ and $(\psi_t)_{t\geq 0}$ on 
  $X$ and $Y$ are \emph{isomorphic} if there is a measure-preserving and essentially invertible
  map $\rho\colon X\to Y$ such that $\psi_t\circ\rho = \rho\circ\varphi_t$ almost everywhere
  for each $t\geq 0$.
  We say that two  Markov lattice semigroups $(T(t))_{t\geq 0}$ and $(S(t))_{t\geq 0}$ on
  $\LL^1(\X)$ and $\LL^1(\Y)$ are \emph{Markov similar} if there 
  is an invertible bi-Markov lattice homomorphism $\Phi\colon\LL^1(\X)\to\LL^1(\Y)$
  such that $S(t)\circ\Phi = \Phi\circ T(t)$ for each $t\geq 0$. The same notions are
  defined for flows and operator groups analogously.
\end{defi}

We call continuous flows (resp.
Koopman semigroups) as described above \emph{topological models} for measurable flows (resp. Markov lattice
semigroups).
See also \cite[Section 12.3]{erg} for this terminology and similar 
results in the time-discrete case. 
For the sake of simplicity, we restrict ourselves to the 
case $p=1$. The idea of the proof of the following result was kindly provided to us by 
Markus Haase. The employed technique of topological models, however, dates back much further.

\begin{thm}[Topological Model]\label{topmod}
Let $A$ be the generator of a strongly continuous Markov lattice semigroup $(T(t))_{t\geq 0}$ on 
a space $\LL^1(\X)$, where $\X = (X, \Sigma, \mu)$ is a finite measure space. Then there exist a compact space $K$, a 
continuous semiflow $(\psi_t)_{t\geq 0}$ on $K$ and a 
%$(\psi_t)_{t\geq 0}$-invariant 
strictly positive Borel probability measure $\nu$ such that the semiflow $(\psi_t)_{t\geq 0}$  
induces a Koopman-semigroup on $\LL^1(K,\nu)$ which is Markov similar to the semigroup 
$(T(t))_{t\geq 0}$ on $\LL^1(\X)$. The measure
$\nu$ is $(\psi_t)_{t\geq 0}$-invariant if and only if $(T(t))_{t\geq 0}$ is a bi-Markov
lattice semigroup.
\end{thm}
\begin{proof}
Consider 
$\mathcal{A}\coloneqq \{ f\in \LL^\infty (\X) \colon s\mapsto T(s)f \text{ is }\norm{\argument}_\infty\text{-continuous}\}$.
Since each operator $T(t)$ is contractive on $\LL^\infty(\X)$ and multiplicative by \Cref{lem:multlatt}, 
$\mathcal{A}$ is an algebra and clearly $\mathds{1}\in \mathcal{A}$.
Furthermore, $\mathcal{A}$ is closed with respect to $\norm{\argument}_\infty$ and closed under conjugation.
Therefore, $\mathcal{A}$ is a commutative $\mathrm{C}^*$-algebra invariant under $(T(t))_{t\geq 0}$.

We show that $\mathcal{A}$ is dense in $\LL^1(\X)$. The strong continuity
of $(T(t))_{t\geq 0}$ on $\LL^1(\X)$ implies that $\norm{\argument}_1$-$\lim_{t \searrow 0} \frac{1}{t}\int_0^t
 T(r)f  \dr = f$ for each $f\in \LL^\infty(\X)$.
Therefore, it suffices to show that $\int_0^t  T(r)f \dr \in \mathcal{A}$ for $f\in \LL^\infty(\X)$.
For all $0\leq s\leq t$ and $f\in \LL^\infty(\X)$ 
\begin{align*}
\abs*{T(s)\int_0^t  T(r)f  \dr-\int_0^t  T(r)f \dr}
&=\abs*{\int_0^t  T(s+r)f \dr - \int_0^t  T(r)f \dr} \\
&=\abs*{\int_s^{t+s}  T(r)f \dr - \int_0^t T(r)f \dr} \\
&\leq \abs*{\int_t^{t+s} T(r)f \dr} +  \abs[\Bigg]{\int_0^s T(r)f \dr} \\
&\leq 2s\norm{f}_{\infty}\mathds{1}
\end{align*}
since each $T(t)$ is $\norm{\argument}_\infty$-contractive.
This shows that $s\mapsto T(s)\int_0^t  T(r)f  \dr$ is continuous at zero and hence on $[0,\infty)$ with respect to 
$\norm{\argument}_\infty$. Therefore,
$\mathcal{A}$ is dense in $\LL^1(\X)$.

By a combination of the Gelfand-Naimark theorem and the Riesz representation theorem 
as in \cite[Section 12.3]{erg} or \cite[Theorem 5.14]{deJeu2017} one obtains a compact space 
$K$, a $*$-isomorphism $\Phi\colon \mathcal{A} \to \mathrm{C}(K)$ with $\Phi\mathds{1}=\mathds{1}$,
a unique probability measure $\nu$ on $K$ such that 
\begin{align*}
  \int_X \Phi^{-1} g \dmu = \int_K g \dnu
\end{align*}
for all $g \in \mathrm{C}(K)$ and a semiflow semiflow $(\psi_t)_{t\geq 0}$ on $K$ such 
that $T(t)|_{\mathcal{A}} = \Phi^{-1}\circ T_{\psi_t} \circ \Phi$.
By \cite[Theorem 4.17]{erg}, the semiflow $\psi$ is continuous, cf. also \cite[Theorem B-II.3.4]{lecture}.

Moreover, $\Phi$ 
extends to a bi-Markov lattice homomorphism $\Phi\colon\mathrm{L}^1(\X) \to \mathrm{L}^1(K, \nu)$.
Let $(S(t))_{t\geq0}$ denote the semigroup $(T(t))_{t\geq0}$ induces on $\mathrm{L}^1(K, \nu)$ 
via $\Phi$. Then 
\begin{align}\label{koopman?}
  S(t)[f]_\nu = [f\circ \psi_t]_\nu
\end{align}
for all continuous functions $f\in\mathrm{C}(K)$. By a standard approximation argument,
this holds for all bounded, Baire-measurable functions, cf. \cite[Theorem E.1]{erg}.
Via monotone approximation, \eqref{koopman?} extends to all positive integrable functions 
and is hence valid for all $[f]_\nu\in\mathrm{L}^1(K, \nu)$.
Finally, $(T(t))_{t\geq 0}$ is bi-Markov if and only if
$(S(t))_{t\geq 0}$ is, if and only if 
\begin{align*}
  \int_K f\dnu = \int_K S(t)f \dnu = \int_K f\circ\psi_t\dnu
\end{align*}
for all $f\in\mathrm{L}^1(K, \nu)$. This is given if and only 
if $(\psi_t)_{t\geq0}$ preserves $\nu$.
\end{proof}

% Recall that a finite measure $\mu$ on $(X,\Sigma)$ is called \emph{separable} if the measure algebra 
% $\Sigma(\X)$ is separable as a metric space with respect to the metric $d(A,B)\coloneqq \mu(A\Delta B)$.
% This is equivalent to the condition that the spaces $\LL^p(\X)$ are separable for all 
% $1\leq p<\infty$, cf.\ \cite[7.14(iv)]{bogachev}.

\begin{prop}\label{metr}
If, in the situation of \Cref{topmod}, the measure space is separable, then the compact space $K$ can be 
chosen to be metrizable.
\end{prop}
\begin{proof}
Let $\mathcal{A}$ be the algebra 
\[ \mathcal{A}\coloneqq \{ f\in \LL^\infty (\X) \colon s\mapsto T(s)f \text{ is 
}\norm{\argument}_\infty\text{-continuous}\}\]
from the proof of \Cref{topmod}. Since $\mathrm{L}^1(\X)$ is separable and 
$\mathcal{A}$ is dense in the former, there is a countable dense subset 
$D_0$ of $\mathcal{A}$. We set 
\begin{align*}
  D \coloneqq \{ T(t)f \colon f\in D_0, t\in\Q_+\} \subset\mathcal{A}
\end{align*}
and denote by $\mathcal{A}_0$ the $\mathrm{C}^*$-subalgebra of 
$\mathcal{A}$ generated by $D$. The algebra $\mathcal{A}_0$ is then
separable, dense in $\mathrm{L}^1(\X)$ and since 
$\mathcal{A}_0 \subset\mathcal{A}$, $T(t)\mathcal{A}_0 \subset \mathcal{A}_0$
for not only $t\in \Q_+$ but $t\in\R_+$. To complete the proof, one can now proceed as in the 
proof of \Cref{topmod} with $\mathcal{A}$ replaced by $\mathcal{A}_0$, obtaining a compact 
representation
space which is metrizable because $\mathrm{C}(K)$ is separable.
\end{proof}

\begin{rem}\label{grouprem}
  With slight notational adjustments, the proofs of the previous two 
  results also work for Markov lattice groups and continuous flows.
\end{rem}

\begin{cor}\label{flowrep}
  Let $\mathrm{X} = (X, \Sigma,\mu)$ be a standard probability space and 
  $(\varphi_t)_{t\in\R}$ a measurable and  measure-preserving flow on $X$.
  Then there are a compact metric space $K$, a continuous flow $(\psi_t)_{t\in\R}$
  on $K$ and
  a strictly positive $(\psi_t)_{t\in\R}$-invariant Borel probability measure $\nu$ on $K$ 
  so that the flows $(\varphi_t)_{t\in\R}$ and $(\psi_t)_{t\in\R}$ are isomorphic.
\end{cor}
\begin{proof}
  By \Cref{measflow} and \Cref{first}, the flow $(\varphi_t)_{t\in\R}$ induces
  a bi-Markov group on $\LL^1(\X)$ and so \Cref{grouprem} shows that there are 
  a compact metric space $K$, a continuous flow 
  $(\psi_t)_{t\in\R}$ and a strictly positive $(\psi_t)_{t\in\R}$-invariant probability measure
  $\nu$ on $K$ such that the groups $(T(t))_{t\in\R}$ and $(S(t))_{t\in\R}$ induced by the 
  flows $(\varphi_t)_{t\in\R}$ and $(\psi_t)_{t\in\R}$ are Markov similar via an invertible bi-Markov
  lattice homomorphism $\Phi$. Applying von Neumann's theorem
  shows that there is a measurable and  measure-preserving
  map $\rho\colon Y\to X$ such that $\Phi = T_\rho$ and $\rho$ is essentially 
  invertible because $\Phi$ is invertible, see \cite[Corollary 7.21]{erg}. The identity
  $\Phi\circ T(t) = S(t)\circ\Phi$ now shows that $\varphi_t\circ\rho = \rho\circ\psi_t$ 
  $\nu$-almost everywhere, see \cite[Proposition 7.19]{erg}.
\end{proof}

\begin{rem}
 \Cref{flowrep} is similar to \cite[Theorem 5]{flows} but for two important differences:
 On the one hand, the authors of \cite{flows} work with a slightly stronger notion 
 of isomorphism of flows. On the other hand, the models considered in \cite{flows} 
 need not be compact.
\end{rem}

\section{Ergodic flows}
\label{sec:ergodicflows}

In this section, we give an operator-theoretic proof for the fact that a measure-preserving 
ergodic flow on a separable measure space comprises at most
countably many non-ergodic maps if it induces a strongly continuous group on 
$\LL^2(\X)$. This is a special case of \cite[Theorem 1]{elements}
where $\R^k$-actions where considered. First, recall the following two properties.

\begin{defi}\label{ergdef}
Let $\X = (X, \Sigma, \mu)$ be a finite measure space.
\begin{enumerate}[(a)]
\item A measure-preserving semiflow $(\varphi_t)_{t\geq 0}$ on $X$ is called 
\emph{ergodic}, if for every $A\in \Sigma$
one has $A\subseteq \varphi_t^{-1}(A)$ modulo null-sets
for each $t\geq 0$ if and only if $A = \emptyset$ or $A = X$ modulo null-sets.
\item A Markov lattice semigroup $(T(t))_{t\geq 0}$ on $\LL^1(\X)$ is called
\emph{irreducible} if there are no nontrivial closed $(T(t))_{t\geq 0}$-invariant ideals in
$\LL^1(\X)$.
\end{enumerate}
\end{defi}

Similar notions can be defined for measure-preserving flows and Markov lattice groups
in a straightforward way and it is not difficult to see that a flow 
$(\varphi_t)_{t\in\R}$ is ergodic if and only if one/both of the associated semiflows 
$(\varphi_t)_{\pm t\geq 0}$ is/are ergodic. As a consequence of 
\cite[Theorem 7.10]{erg}, the two properties in \Cref{ergdef} are equivalent 
for
a (semi)flow and its corresponding Koopman (semi)group.
Obviously, if some $\varphi_{s}$ is ergodic (meaning that $A\subseteq \varphi_{s}^{-1}(A)$ 
implies $\mu(A)\in\{0,1\}$), so is the semiflow $(\varphi_t)_{t\geq0}$, while the converse is not 
true in general. However, it follows from \cite[Theorem 1]{elements} that for an ergodic \emph{flow} 
$(\varphi_t)_{t\in\R}$ all but at most countably many maps $\varphi_t$ are ergodic. In the 
following we give a short proof of this fact using the Perron-Frobenius spectral theory.
We shall prove this for semiflows for which the induced Koopman semigroup is strongly continuous 
on $\LL^2(\X)$. In light of the previous remarks, this implies 
\cite[Theorem 1]{elements}.

\begin{thm}\label{erg}
Let $(\varphi_t)_{t\geq 0}$ be a measure-preserving semiflow on a separable finite measure 
space $\X=(X,\Sigma,\mu)$ such that the corresponding Koopman semigroup $(T(t))_{t\geq 0}$ is strongly continuous on 
$\LL^2(\X)$. 
If $(\varphi_t)_{t\geq 0}$ is ergodic, then at most countably many maps $\varphi_{t}$ 
are not ergodic. 
\end{thm}
\begin{proof}
Let $(\varphi_t)_{t\geq 0}$ be measure-preserving and ergodic and its induced Koopman semigroup $(T(t))_{t\geq 0}$ strongly 
continuous
on $\LL^2(\X)$ with generator $A$. Since the eigenspaces of an isometry on $\LL^2(\X)$ are 
pairwise orthogonal, the point-spectrum of $T(t)$ is countable for each $t \geq 0$ because $\LL^2(\X)$ is
separable. The spectral inclusion theorem \cite[Theorem IV.3.6]{engnag} implies that 
the boundary point spectrum $G\coloneqq \mathrm{P}\sigma(A) \cap \mathrm{i}\R$ is countable.

For a fixed $s > 0$ the map  $\varphi_s$ is not ergodic if and only if 
$\dim  \fix(T(s))\geq 2$, see \cite[Proposition 7.15]{erg}. By the spectral mapping theorem for the point spectrum
\cite[Corollary IV.3.8]{engnag}
\[
\fix(T(s))=\overline{\lin}\{f\in \LL^2(\X) \colon Af=\mathrm{i}\lambda f \text{ for some }\mathrm{i}\lambda \in G \text{ 
with 
}\mathrm{e}^{\mathrm{i}\lambda s}=1\}.
\]
Moreover, for an irreducible Markov lattice semigroup the common fixed space 
\begin{align*}
 \bigcap_{t\geq 0} \fix(T(t))=\ker A
\end{align*}
is one-dimensional since it is a Banach sublattice of $\LL^2(\X)$ and for real-valued $f\in\bigcap_{t\geq 0} \fix(T(t))$
all characteristic functions $\1_{[f > c]}$, $c\in\R$, are contained in the fixed space (use \cite[Exercise 7.13]{erg}).
Thus, $\dim \fix(T(s)) \geq 2$ if and only if there exists $\mathrm{i}\lambda \in G$ such 
that $\lambda s = 2\pi k$ for some $k\in\mathds{Z}\setminus\{0\}$.
Since $G$ is countable, this can only be the case for at most countably many $s>0$.
\end{proof}

As a consequence of \Cref{erg}, \Cref{measflow} and the remarks made above, we note the following corollary.

\begin{cor}
Let $(\varphi_t)_{t\in\R}$ be a measurable and measure-preserving flow on a separable finite measure 
space $\X$.
If $(\varphi_t)_{t\in\R}$ is ergodic, then at most countably many maps $\varphi_{t}$ 
are not ergodic. 
\end{cor}

 The assertion in \Cref{erg} does not hold for non-separable spaces. 
 Take for example the Bohr compactification $\mathrm{b}\R$ of the additive group of real numbers and define
 for $g\in \mathrm{b}\R$ the translation
 $\varphi_g\colon \mathrm{b}\R\to \mathrm{b}\R$, $h \mapsto h + g$. Then $(\varphi_t)_{t\in\R}$ is a continuous
 flow on $\mathrm{b}\R$ preserving the Haar measure $\mathrm{m}$. Denote the corresponding Koopman 
 group on $\LL^2(\mathrm{b}\R, \mathrm{m})$ by $(T(t))_{t\in\R}$. If $f\in\LL^2(\mathrm{b}\R, \mathrm{m})$
 is such that $T(t)f = f\circ\varphi_t = f$ for all $t\in\R$, then $f\circ\varphi_g = f$ for all $g\in\mathrm{b}\R$ since
 $\R$ is dense in $\mathrm{b}\R$. It follows that $f$ is constant almost everywhere.
 Because $T(t)f = f$ is equivalent to $T(-t)f = f$ for each $t\geq 0$, this shows 
 that the flow $(\varphi_t)_{t\geq 0}$ is ergodic. However, $\varphi_t$ is not ergodic 
 since $T(t)f = f$ for all periodic functions with period $t$.

\subsection*{Acknowledgement}We express our sincere gratitude towards the anonymous referee for 
his insightful observations and detailed comments.
 
\printbibliography

\end{document}